\documentclass[final]{siamltex}
\DeclareOldFontCommand{\sc}{\normalfont\scshape}{\@nomath\sc}
\usepackage{amsmath}
\usepackage{amssymb}
\usepackage{array}
\usepackage{blkarray}
\usepackage{caption}
\usepackage{capt-of}
\usepackage{color}
\usepackage{colortbl}
\usepackage{fancyhdr}
\usepackage{fancyvrb}
\usepackage{float}
\usepackage[T1]{fontenc}
\usepackage{framed}
\usepackage{geometry}
\usepackage{graphicx}
\usepackage[latin1]{inputenc}
\usepackage{stmaryrd}
\usepackage{upquote}
\usepackage{url}
\usepackage[usenames,dvipsnames,svgnames,table]{xcolor}
\usepackage{xspace}	
\usepackage{booktabs}
\usepackage{tikz}
\usepackage{mathdots}
\usepackage{soul}
\definecolor{lightgray}{gray}{0.5}
\definecolor{gray}{gray}{0.7}
\definecolor{darkgreen}{rgb}{0,0.6,0.13}
\newcommand{\nc}{\newcommand}
\nc{\dsp}{\displaystyle}
\nc{\txt}{\textstyle}
\nc{\reff}[1]{(\ref{#1})}
\nc{\mrm}[1]{\mathrm{#1}}
\nc{\udl}[1]{\underline{#1}}
\nc{\ovl}[1]{\overline{#1}}
\nc{\al}{\underline{\boldsymbol{\alpha}}}
\nc{\la}{\underline{\boldsymbol{\lambda}}}
\nc{\llbr}{\llbracket}
\nc{\rrbr}{\rrbracket}
\nc{\lbr}{\lbrack}
\nc{\rbr}{\rbrack}
\nc{\N}{\mathbb{N}}
\nc{\Z}{\mathbb{Z}}
\nc{\D}{\mathbb{D}}
\nc{\Q}{\mathbb{Q}}
\nc{\R}{\mathbb{R}}
\nc{\C}{\mathbb{C}}
\nc{\T}{\mathbb{T}}
\nc{\Stwo}{\mathbb{S}^2}
\nc{\tld}[1]{\tilde{#1}}
\nc{\wtld}[1]{\widetilde{#1}}
\nc{\hu}{\hat{u}}
\nc{\wh}[1]{\widehat{#1}}
\nc{\Fbf}{\textbf{F}}
\nc{\Gbf}{\textbf{G}}
\nc{\Lbf}{\textbf{L}}
\nc{\Nbf}{\textbf{N}}
\nc{\Ibf}{\textbf{I}}
\nc{\Dbf}{\textbf{D}} 
\nc{\Tbf}{\textbf{T}}
\nc{\Jbf}{\textbf{J}} 
\nc{\Rbf}{\textbf{R}}   
\nc{\ph}{\varphi}
\nc{\sumeven}{\sum_{k=-N/2}^{N/2}{\hspace{-0.3cm}}'{\;\,}}
\nc{\sumevenk}{\sum_{k=-n/2}^{n/2}{\hspace{-0.3cm}}'{\;\,}}
\nc{\sumevenj}{\sum_{j=-m/2}^{m/2}{\hspace{-0.3cm}}'{\;\,}}
\nc{\sumodd}{\sum_{k=-\frac{N-1}{2}}^{\frac{N-1}{2}}}
\nc{\sumoddl}{\sum_{l=-\frac{N-1}{2}}^{\frac{N-1}{2}}}
\nc{\cqfd}{~\hbox{\vrule width 2.5pt depth 2.5 pt height 3.5 pt}}
\nc{\bs}[1]{\boldsymbol{#1}}
\nc{\rr}[1]{\textcolor{red}{#1}}
\nc{\bb}[1]{\textcolor{blue}{\xout{#1}}}
\nc{\ra}[1]{\renewcommand{\arraystretch}{#1}}

\setulcolor{red}
\soulregister\cite7
\soulregister\ref7
\soulregister\pageref7
\title{New error bounds for deep ReLU networks using sparse grids}
\author{Hadrien Montanelli\thanks{Department of Applied Physics and Applied Mathematics, Columbia University, New York,
NY 10027, United States.\ \ This research is supported in part by US NSF DMS-1719699 and the NSF TRIPODS program CCF-1704833.}
\and Qiang Du\thanks{Department of Applied Physics and Applied Mathematics, Columbia University, New York,
NY 10027, United States.\ \ This research is supported in part by US NSF DMS-1719699 and the NSF TRIPODS program CCF-1704833.}
}

\begin{document}

\maketitle

%%%%%%%%%%%%%%%%%%%%%%%%%%%%%%%%%%%%%%%%%%%%%%%%%%%%%%%%%%%%%%%%%%%%%%%%%%%
\begin{abstract}
We prove a theorem concerning the approximation of multivariate functions by deep ReLU networks.
We present new error estimates for which the curse of the dimensionality is lessened by establishing a connection
with sparse grids.
\end{abstract}

%%%%%%%%%%%%%%%%%%%%%%%%%%%%%%%%%%%%%%%%%%%%%%%%%%%%%%%%%%%%%%%%%%%%%%%%%%%

%%%%%%%%%%%%%%%%%%%%%%%%%%%%%%%%%%%%%%%%%%%%%%%%%%%%%%%%%%%%%%%%%%%%%%%%%%%

%%%%%%%%%%%%%%%%%%%%%%%%%%%%%%%%%%%%%%%%%%%%%%%%%%%%%%%%%%%%%%%%%%%%%%%%%%%
\pagestyle{myheadings}
\thispagestyle{plain}

\markboth{MONTANELLI AND DU}{NEW ERROR BOUNDS FOR DEEP RELU NETWORKS}

%%%%%%%%%%%%%%%%%%%%%%%%%%%%%%%%%%%%%%%%%%%%%%%%%%%%%%%%%%%%%%%%%%%%%%%%%%%
\section{Introduction}

Deep learning has been successfully applied to many fields including computer vision, 
speech recognition, and natural language processing \cite{lecun2015}.
It is based on approximations by \textit{deep} networks, as opposed to \textit{shallow} networks.
The latter are neural networks with a single \textit{layer} and correspond to approximations $f_N$ of multivariate 
functions $f:\R^d\rightarrow\R$ of the form
\begin{equation}
f_N(\boldsymbol{x}) = \sum_{i=1}^N \alpha_i \sigma(\boldsymbol{w}_i^T\boldsymbol{x} + \theta_i), 
\quad \alpha_i,\theta_i\in\R, \, \boldsymbol{x},\boldsymbol{w}_i\in\R^d,
\label{shallow}
\end{equation}

\noindent for some \textit{activation function} $\sigma:\R\rightarrow\R$.
The former are neural networks with one or more layers, where each unit of each layer 
performs an operation of the form $\sigma(\boldsymbol{w}\cdot \boldsymbol{x} + \theta)$.
The \textit{depth} of a network is the number of layers, and the \textit{size} is the total number of units.
Shallow networks have depth $1$ and their size is the number $N$ in \eqref{shallow},
while deep networks usually have depth $\gg 1$.
Deep \textit{ReLU} networks use the activation function $\sigma(x) = \max(0,x)$.

One of the most important theoretical problems is to determine why and when deep (but not shallow) networks
can lessen or break the ``curse of dimensionality,'' expression first coined by Bellman in \cite{bellman1961}.
A possible way of addressing this problem is to focus on a particular set of functions which have a very
special structure (such as compositional or polynomial), and to show that for this particular set deep networks perform
extremely well \cite{cohen2016, eldan2016, liang2016, poggio2017, telgarsky2016}.
We follow a different route. 
We consider a space of functions that is more generic for multivariate approximation in high dimensions, 
and prove new error estimates for which the curse of dimensionality is lessened by establishing a connection with 
\textit{sparse grids} \cite{bungartz2004}. 

The theory of approximating functions using neural networks goes back to the late 1980s with the first 
\textit{density results}~\cite{cybenko1989, hornik1989}.
Cybenko \cite{cybenko1989} showed that any continuous functions can be approximated by shallow networks,
while Hornik et al.~\cite{hornik1989} proved a similar result for Borel measurable functions.\footnote{These results automatically 
apply to deep networks since deep networks include networks with a singe layer.}
They both used \textit{sigmoid} functions as activation functions.\footnote{A function $\sigma:\R\rightarrow[0,1]$ 
is said to be a \textit{sigmoid} function if it is non-decreasing,
$\lim_{x\rightarrow-\infty}\sigma(x)=0$ and $\lim_{x\rightarrow+\infty}\sigma(x)=1$, e.g., $\sigma(x)=1/(1+e^{-x})$.
This does not include the ReLU activation function.
Some density results for ReLU functions can be found in the review of Pinkus (e.g., \cite[Prop.~3.7]{pinkus1999}).}

\begin{table}[t]
\caption{\textit{Approximation results with curse of dimensionality.
The first line lists function spaces, the second line the depth and the size
of the network to get accuracy $\epsilon$, 
and the third line the norm used for the estimates, together with our favorite reference.
The domain is either an hypercube, the unit ball $B^d$ of $\R^d$ with respect to the two-norm of vectors, or
a poly-ellipse $E_\rho$ with parameter $\rho$, a generalization of the one-dimensional Bernstein ellipse.}}
\centering
\ra{1.3}
\begin{tabular}{ccc}
\toprule
& \textbf{Shallow} & \textbf{Deep} \\
\midrule
$\boldsymbol{\sigma\in C^\infty(\R)}$ & $f\in W^{m,p}([-1,1]^d)$ & \\
\textbf{(not polynomial)} & depth $1$, size $\mathcal{O}(\epsilon^{-\frac{d}{m}})$ & \textbf{--} \\
& $\Vert\cdot\Vert_p$, \cite[Th.~2.1]{mhaskar1996} & \\\\
$\boldsymbol{\sigma\in C^\infty(\R)}$  & $f$ analytic in $E_\rho$ & \\
\textbf{(not polynomial)} & depth $1$, size $\mathcal{O}(\vert\log_\rho\epsilon\vert^d)$ & \textbf{--} \\
& $\Vert\cdot\Vert_p$, \cite[Th.~2.3]{mhaskar1996} & \\\\
$\boldsymbol{\sigma}$ \textbf{ReLU} & $f\in W^{m,2}(B^d)$ & $f\in W^{m,\infty}([0,1]^d)$ \\
& depth $1$, size $\mathcal{O}(\epsilon^{-\frac{d}{m}})$ & depth $\mathcal{O}(\vert\log_2\epsilon\vert)$, 
size $\mathcal{O}(\epsilon^{-\frac{d}{m}}\vert\log_2\epsilon\vert)$\\
& $\Vert\cdot\Vert_2$, \cite[Cor.~6.10]{pinkus1999} & $\Vert\cdot\Vert_\infty$, \cite[Th.~1]{yarotsky2017} \\
\bottomrule
\end{tabular}
\label{table:curse}
\end{table}

These results were important from a theoretical point of view but what ultimately matters in practice 
is \textit{how fast} approximations by neural networks converge.
For example, for a real-valued function $f$ in $\R^d$ whose smoothness is characterized by some integer~$m$ (typically the order of integrable or bounded derivatives), 
{and for some prescribed accuracy $\epsilon>0$, one tries to show that there exists a neural network $f_N$ of size $N$ that satisfies
\begin{equation}
\Vert f - f_N \Vert \leq \epsilon \quad \text{with} \quad N=\mathcal{O}(\epsilon^{-\frac{d}{m}}),
\footnote{We recall that $N=\mathcal{O}(\epsilon^{-\frac{d}{m}})$ means that there exists a constant $C>0$ such that
$N\leq C\epsilon^{-\frac{d}{m}}$ for sufficiently small values of $\epsilon$. 
As is the case with sparse grids, we will abusively use this notation even if $C=C(d)$ depends on the dimension $d$.}
\label{approx}
\end{equation}

\noindent for some norm $\Vert\cdot\Vert$.
For deep networks, one also wants to find the asymptotic behavior of the depth as a function of the accuracy $\epsilon$.
Results of the form~\eqref{approx} are standard approximation results that suffer from the curse of dimensionality.
For a small dimensions $d$, the size $N$ of the network increases at a reasonable rate as $\epsilon$ goes to zero.
However, $N$ grows \textit{geometrically} with $d$.

Many results of the form~\eqref{approx} have been derived for shallow and deep networks; we list some of them in Table~\ref{table:curse},
for different activation functions $\sigma$.
The function spaces considered in the table are the Sobolev spaces $W^{m,p}(\Omega)$ (for some compact subset $\Omega\subset\R^d$) 
of functions that have (weak) partial derivatives in $L^p(\Omega)$ up to order $m$, i.e.,
\begin{equation}
W^{m,p}(\Omega) = \big\{ f \in L^p(\Omega) : D^{\boldsymbol{k}} f \in L^p(\Omega), \, \vert \boldsymbol{k}\vert_1\leq m \big\},
\label{Sobolev}
\end{equation}

\noindent with multi-index $\boldsymbol{k}=(k_1,\ldots,k_d)\in\N^d$, $\vert \boldsymbol{k} \vert_1 = \sum_{j=1}^d k_j$, and partial derivatives
\begin{equation}
D^{\boldsymbol{k}} f = \frac{\partial^{\vert \boldsymbol{k}\vert_1} f}{\partial x_1^{k_1} \ldots \partial x_d^{k_d}}.
\end{equation}

\noindent These are Banach spaces too corresponding to the completion of $C^m(\Omega)$ (functions
with continuous partial derivatives up to order $m$) with respect to the norm
\begin{equation}
\Vert f \Vert_{m,p} = \left\{
\begin{array}{ll}
\dsp \Bigg(\sum_{0\leq \vert \boldsymbol{k} \vert_1\leq m} \Vert D^{\boldsymbol{k}} f\Vert_p^p\Bigg)^{\frac{1}{p}}, & 1\leq p < \infty, \\[16pt]
\dsp \underset{0\leq\vert \boldsymbol{k} \vert_1\leq m}{\max} \Vert D^{\boldsymbol{k}} f\Vert_\infty, & p = \infty,
\end{array}
\right.
\end{equation}

\noindent with standard $L^p(\Omega)$-norm $\Vert\cdot\Vert_p$.

\begin{table}[t]
\caption{\textit{Approximation results without curse of dimensionality.
The first line lists the function spaces, the second line the depth and the size
of the network to get accuracy $\epsilon$, 
and the third line the norm used for the estimates, together with our favorite reference.
Deep networks can exploit the compositional structure of functions to break the curse of dimensionality;
shallow networks cannot.}}
\centering
\ra{1.3}
\begin{tabular}{ccc}
\toprule
& \textbf{Shallow} & \textbf{Deep} \\
\midrule
$\boldsymbol{\sigma\in C^\infty(\R)}$ & $f\in W^{m,\infty}([-1,1]^d)$, compositional & $f\in W^{m,\infty}([-1,1]^d)$, compositional \\
\textbf{(not polynomial)} & depth $1$, size $\mathcal{O}(\epsilon^{-\frac{d}{m}})$ & depth $\log_2 d$, size $\mathcal{O}((d-1)\epsilon^{-\frac{2}{m}})$ \\
& $\Vert\cdot\Vert_\infty$, \cite[Th.~1]{poggio2017} & $\Vert\cdot\Vert_\infty$, \cite[Th.~2]{poggio2017} \\\\
$\boldsymbol{\sigma}$ \textbf{ReLU} & $f$ Lipschitz, $[-1,1]^d$, compositional & $f$ Lipschitz, $[-1,1]^d$, compositional \\
& depth $1$, size $\mathcal{O}(\epsilon^{-d})$ & depth $\log_2 d$, size $\mathcal{O}((d-1)\epsilon^{-2})$ \\
& $\Vert\cdot\Vert_\infty$, \cite[Th.~4]{poggio2017} & $\Vert\cdot\Vert_\infty$, \cite[Th.~4]{poggio2017} \\
& & (see also \cite[Cor.~10]{liang2016}) \\
\bottomrule
\end{tabular}
\label{table:nocurse}
\end{table}

As we mentioned previously, some results without curse of dimensionality for deep (but not shallow) networks have been recently derived; 
some of them are listed in Table~\ref{table:nocurse}.
To derive such results, one has to consider functions with a very special structure.
For example, in \cite{poggio2017}, the authors deal with \textit{compositional} functions.
This class includes in particular functions that are a composition of two-dimensional functions, e.g., in dimension $d=8$,
\begin{equation}
f(x_1,\ldots,x_8) = f_3\Big(f_{21}\big(f_{11}(x_1,x_2), f_{12}(x_3,x_4)\big), f_{22}\big(f_{13}(x_5,x_6),f_{14}(x_7,x_8)\big)\Big),
\end{equation}

\noindent for some bivariate functions $f_{11}$, $f_{12}$, $f_{13}$, $f_{14}$, $f_{21}$, $f_{22}$ and $f_3$. 
Such functions can be represented by a binary tree with $d$ inputs variables, $\log_2 d$ levels
and $(d-1)$ nodes, and each of the nodes can be approximated by a subnetowrk of size $N/(d-1)$; see Figure~\ref{compositional}.

\begin{figure}[t]
\centering
\begin{tikzpicture}
[   cnode/.style={draw=black, fill=#1, minimum width=3mm, circle},
]
    \node[cnode = black, label = 270: $x_1$] (x1) at (1, -3) {};
    \node[cnode = black, label = 270: $x_2$] (x2) at (2, -3) {};
    \node[cnode = black, label = 270: $x_3$] (x3) at (3, -3) {};
    \node[cnode = black, label = 270: $x_4$] (x4) at (4, -3) {};
    \node[cnode = black, label = 270: $x_5$] (x5) at (5, -3) {};
    \node[cnode = black, label = 270: $x_6$] (x6) at (6, -3) {};
    \node[cnode = black, label = 270: $x_7$] (x7) at (7, -3) {};
    \node[cnode = black, label = 270: $x_8$] (x8) at (8, -3) {};
    \node[cnode = gray, label = 270: $f_{11}$] (y1) at (1.5, -2) {};
    \node[cnode = gray, label = 270: $f_{12}$] (y2) at (3.5, -2) {};
    \node[cnode = gray, label = 270: $f_{13}$] (y3) at (5.5, -2) {};
    \node[cnode = gray, label = 270: $f_{14}$] (y4) at (7.5, -2) {};
    \node[cnode = gray, label = 270: $f_{21}$] (z1) at (2.5, -1) {};
    \node[cnode = gray, label = 270: $f_{22}$] (z2) at (6.5, -1) {};
    \node[cnode = gray, label = 270: $f_3$] (f) at (4.5, 0) {};
    \draw (x1) -- (y1) node [below, sloped, pos=0.5] {};
    \draw (x2) -- (y1) node [below, sloped, pos=0.5] {};
    \draw (x3) -- (y2) node [below, sloped, pos=0.5] {};
    \draw (x4) -- (y2) node [below, sloped, pos=0.5] {};
    \draw (x5) -- (y3) node [below, sloped, pos=0.5] {};
    \draw (x6) -- (y3) node [below, sloped, pos=0.5] {};
    \draw (x7) -- (y4) node [below, sloped, pos=0.5] {};
    \draw (x8) -- (y4) node [below, sloped, pos=0.5] {};
    \draw (y1) -- (z1) node [below, sloped, pos=0.5] {};
    \draw (y2) -- (z1) node [below, sloped, pos=0.5] {};
    \draw (y3) -- (z2) node [below, sloped, pos=0.5] {};
    \draw (y4) -- (z2) node [below, sloped, pos=0.5] {};
    \draw (z1) -- (f) node [below, sloped, pos=0.5] {};
    \draw (z2) -- (f) node [below, sloped, pos=0.5] {};
    \node[draw=none] (ghost1) at (9, -2.25) {};
    \node[draw=none] (ghost2) at (9, 0.25) {};
    \draw [thick, <->] (ghost1) -- node[right]{$\log_2 d=3$} (ghost2);
\end{tikzpicture}
\caption{\textit{The binary tree of a compositional function in dimension $d=8$. The tree has $\log_2 d=3$ levels and $d-1=7$ nodes/local functions.
To approximate such a function with a deep network of size $N$, one can mimic the compositional structure by using a network with depth $\log_2 d$
and  made of $7$ subnetworks, each subnetwork approximating a local function with $N/7$ units.}}
\label{compositional}
\end{figure}
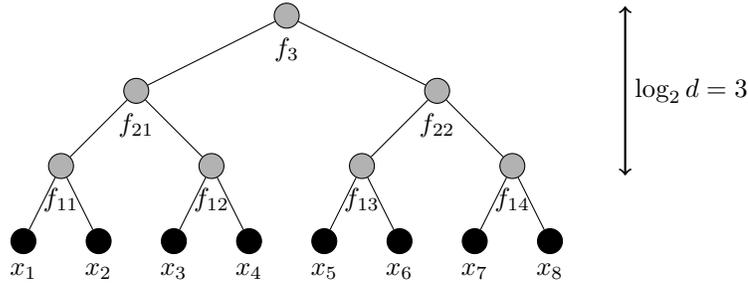

We will show in this paper that functions in the so-called Korobov spaces $X^{2,p}(\Omega)$ (for some compact subset $\Omega\subset\R^d$) 
of \textit{mixed derivatives} of order two can be represented to accuracy $\epsilon$ by deep networks of depth $\mathcal{O}(\vert\log_2\epsilon\vert\log_2d)$ 
and size
\begin{equation}
N=\mathcal{O}(\epsilon^{-\frac{1}{2}}\vert\log_2\epsilon\vert^{\frac{3}{2}(d-1)+1}\log_2d).
\label{size}
\end{equation}

\noindent The estimate is achieved via sparse grid approximations to functions in $X^{2,p}(\Omega)$. 
The curse of dimensionality is not totally overcome but is significantly lessened since the exponent $d$ only affects
logarithmic factors $\vert\log_2\epsilon\vert$.\footnote{Let us emphasize that the constant in \eqref{size}, 
however, still depends exponentially on the dimension $d$.}
As we will see in Section \ref{section2.2}, Korobov spaces $X^{2,p}(\Omega)$ are subsets of Sobolev spaces $W^{2,p}(\Omega)$.

The reminder of the paper is structured as follows.
We review Korobov spaces and sparse grids in Section \ref{section2}, and prove our theorem in Section \ref{section3}.

%%%%%%%%%%%%%%%%%%%%%%%%%%%%%%%%%%%%%%%%%%%%%%%%%%%%%%%%%%%%%%%%%%%%%%%%%%%
\section{Korobov spaces and sparse grids}\label{section2}

We review in this section Korobov spaces and sparse grids, which go back to Korobov \cite{korobov1959} and Smolyak \cite{smolyak1963}, 
and were rediscovered by Zenger in \cite{zenger1991} for solving partial differential equations. 
Since 1990, Korobov spaces/sparse grids and related \textit{hyperbolic cross} approximation \cite{schmeisser2004, shen2010}
have been used extensively in the context of high-dimensional function approximation, and the group of Michael Griebel has been particularly influential.
For details, we recommend the exhaustive review \cite{bungartz2004}.

%%%%%%%%%%%%%%%%%%%%%%%%%%%%%%%%%%%%%%%%%%%%%%%%%%%%%%%%%%%%%%%%%%%%%%%%%%%
\subsection{One-dimensional hierarchical basis}\label{section2.1}

The first ingredient for cooking sparse grids is a \textit{hierarchical basis} of functions.
To approximate functions of one variable $x$ on $\Omega = [0, 1]$, one considers a family of grids $\Omega_l$ 
of \textit{level}~$l$ characterized by a grid size $h_l=2^{-l}$ and $2^l-1$ points $x_{l,i}=i h_l$, $1\leq i\leq 2^l-1$.
For each $\Omega_l$, ones considers piecewise linear hat functions $\phi_{l,i}$ centered at $x_{l,i}$ defined by
\begin{equation}
\phi_{l,i}(x) = \phi\bigg(\frac{x - x_{l,i}}{h_l}\bigg), \quad 1\leq i\leq 2^l-1,
\label{piecewiselin}
\end{equation}

\noindent where $\phi$ is the mother of all hat functions,
\begin{equation}
\phi(x) = \left\{
\begin{array}{ll}
1 - \vert x\vert, & \text{if}\,\,x\in[-1,1], \\[16pt]
0, & \text{otherwise}.
\end{array}
\right.
\end{equation}

\noindent Note that $\Vert\phi_{l,i}\Vert_\infty\leq1$ for all $l$ and $i$. 
One then considers the function spaces $V_l$ spanned by such functions,
\begin{equation}
V_l = \mrm{span}\{\phi_{l,i} : 1\leq i\leq 2^{l}-1\},
\end{equation}

\noindent and the \textit{hierarchical increment} spaces $W_l$ given by
\begin{equation}
W_l = \mrm{span}\{\phi_{l,i} : i\in I_l\},
\end{equation}

\noindent where $I_l = \{i\in\N : 1\leq i\leq 2^{l}-1, \, i\,\text{odd}\}$.
These increment spaces satisfy the relation
\begin{equation}
V_n = \bigoplus_{1\leq l\leq n} W_l.
\end{equation}

\noindent The basis that corresponds to the $W_l$'s for $1\leq l\leq n$ is called the \textit{hierarchical basis},
while the basis of $V_n$ is called the \textit{nodal basis}. We show both bases for $n=3$ in Figure~\ref{hierarchical}.

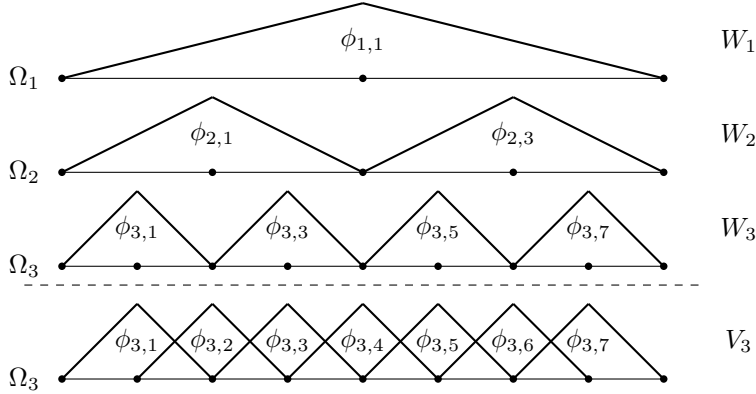
\begin{figure}[t]
\centering
\begin{tikzpicture}
[   cnode/.style={draw=black, fill=#1, minimum width=3mm, circle},
]   
	   % W1:
      \draw (0, 0) -- (8, 0) {};
      \draw[thick] (0, 0) -- (4, 1) {};
      \draw[thick] (4, 1) -- (8, 0) {};
      \node[fill, circle, scale=0.3] at (0, 0) {};
      \node[fill, circle, scale=0.3] at (4, 0) {};
      \node[fill, circle, scale=0.3] at (8, 0) {};
      \node at (-0.5, 0) {$\Omega_1$};
      \node at (9, 0.5) {$W_1$};
      \node at (4, 0.5) {$\phi_{1,1}$};
       
      % W2:
      \draw (0, -1.25) -- (8, -1.25) {};
      \draw[thick] (0, -1.25) -- (2, -0.25) {};
      \draw[thick] (2, -0.25) -- (4, -1.25) {};
      \draw[thick] (4, -1.25) -- (6, -0.25) {};
      \draw[thick] (6, -0.25) -- (8, -1.25) {};
      \node[fill, circle, scale=0.3] at (0, -1.25) {};
      \node[fill, circle, scale=0.3] at (2, -1.25) {};
      \node[fill, circle, scale=0.3] at (4, -1.25) {};
      \node[fill, circle, scale=0.3] at (6, -1.25) {};
      \node[fill, circle, scale=0.3] at (8, -1.25) {};
      \node at (-0.5, -1.25) {$\Omega_2$};
      \node at (9, -0.75) {$W_2$};
      \node at (2, -0.75) {$\phi_{2,1}$};
      \node at (6, -0.75) {$\phi_{2,3}$};
      
      % W3:
      \draw (0, -2.5) -- (8, -2.5) {};
      \draw[thick] (0, -2.5) -- (1, -1.5) {};
      \draw[thick] (1, -1.5) -- (2, -2.5) {};
      \draw[thick] (2, -2.5) -- (3, -1.5) {};
      \draw[thick] (3, -1.5) -- (4, -2.5) {};
      \draw[thick] (4, -2.5) -- (5, -1.5) {};
      \draw[thick] (5, -1.5) -- (6, -2.5) {};
      \draw[thick] (6, -2.5) -- (7, -1.5) {};
      \draw[thick] (7, -1.5) -- (8, -2.5) {};
      \node[fill, circle, scale=0.3] at (0, -2.5) {};
      \node[fill, circle, scale=0.3] at (1, -2.5) {};
      \node[fill, circle, scale=0.3] at (2, -2.5) {};
      \node[fill, circle, scale=0.3] at (3, -2.5) {};
      \node[fill, circle, scale=0.3] at (4, -2.5) {};
      \node[fill, circle, scale=0.3] at (5, -2.5) {};
      \node[fill, circle, scale=0.3] at (6, -2.5) {};
      \node[fill, circle, scale=0.3] at (7, -2.5) {};
      \node[fill, circle, scale=0.3] at (8, -2.5) {};
      \node at (-0.5, -2.5) {$\Omega_3$};
      \node at (9, -2) {$W_3$};
      \node at (1, -2) {$\phi_{3,1}$};
      \node at (3, -2) {$\phi_{3,3}$};
      \node at (5, -2) {$\phi_{3,5}$};
      \node at (7, -2) {$\phi_{3,7}$};
     
       % V3:
      \draw[dashed] (-0.5, -2.75) -- (8.5, -2.75) {};
      \draw (0, -4) -- (8, -4) {};
      \draw[thick] (0, -4) -- (1, -3) {};
      \draw[thick] (1, -3) -- (2, -4) {};
      \draw[thick] (1, -4) -- (2, -3) {};
      \draw[thick] (2, -3) -- (3, -4) {};
      \draw[thick] (2, -4) -- (3, -3) {};
      \draw[thick] (3, -3) -- (4, -4) {};
      \draw[thick] (3, -4) -- (4, -3) {};
      \draw[thick] (4, -3) -- (5, -4) {};
      \draw[thick] (4, -4) -- (5, -3) {};
      \draw[thick] (5, -3) -- (6, -4) {};
      \draw[thick] (5, -4) -- (6, -3) {};
      \draw[thick] (6, -3) -- (7, -4) {};
      \draw[thick] (6, -4) -- (7, -3) {};
      \draw[thick] (7, -3) -- (8, -4) {};
      \node[fill, circle, scale=0.3] at (0, -4) {};
      \node[fill, circle, scale=0.3] at (1, -4) {};
      \node[fill, circle, scale=0.3] at (2, -4) {};
      \node[fill, circle, scale=0.3] at (3, -4) {};
      \node[fill, circle, scale=0.3] at (4, -4) {};
      \node[fill, circle, scale=0.3] at (5, -4) {};
      \node[fill, circle, scale=0.3] at (6, -4) {};
      \node[fill, circle, scale=0.3] at (7, -4) {};
      \node[fill, circle, scale=0.3] at (8, -4) {};
      \node at (-0.5, -4) {$\Omega_3$};
      \node at (9, -3.5) {$V_3$};
      \node at (1, -3.5) {$\phi_{3,1}$};
      \node at (2, -3.5) {$\phi_{3,2}$};
      \node at (3, -3.5) {$\phi_{3,3}$};
      \node at (4, -3.5) {$\phi_{3,4}$};
      \node at (5, -3.5) {$\phi_{3,5}$};
      \node at (6, -3.5) {$\phi_{3,6}$};
      \node at (7, -3.5) {$\phi_{3,7}$};
       
\end{tikzpicture}
\caption{\textit{Piecewise linear hierarchical basis (top) and nodal basis (bottom) for $n=3$.
The hierarchical basis consists of functions $\phi_{1,1}$, $\phi_{2,1}$, $\phi_{2,3}$, $\phi_{3,1}$, $\phi_{3,3}$, $\phi_{3,5}$ and $\phi_{3,7}$,
which live on three different grids.
Note that for each grid the supports of the hierarchical basis functions are mutually disjoint.
The nodal basis consist of functions $\phi_{3,i}$, $1\leq i\leq 7$, which live on the same grid.
Both bases span the same space of dimension $2^n-1=7$ since $V_3 = W_1 \bigoplus W_2 \bigoplus W_3$.}}
\label{hierarchical}
\end{figure}

Let us conclude this subsection by mentioning that one-dimensional hierarchical bases are not limited to piecewise linear functions \eqref{piecewiselin}.
These can be generalized to piecewise higher-order polynomials \cite[Th.~4.8]{bungartz2004},
and also to other multiscale bases such as wavelets \cite{grossmann1984} (see, e.g., \cite{sprengel1998}).

%%%%%%%%%%%%%%%%%%%%%%%%%%%%%%%%%%%%%%%%%%%%%%%%%%%%%%%%%%%%%%%%%%%%%%%%%%%
\subsection{Multi-dimensional hierarchical basis and Korobov spaces}\label{section2.2}

The second ingredient is to employ a \textit{tensor product} construction to approximate functions of $d$ variables 
$\boldsymbol{x}=(x_1,\ldots,x_d)$ in $\Omega = [0, 1]^d$.
One considers a family of grids $\Omega_{\boldsymbol{l}}$ of level $\boldsymbol{l}=(l_1,\ldots,l_d)$ with points 
$\boldsymbol{x}_{\boldsymbol{l},\boldsymbol{i}} = \boldsymbol{i}\cdot \boldsymbol{h_l}$, 
$\boldsymbol{1} \leq \boldsymbol{i} \leq \boldsymbol{2^l-1}$, obtained by the tensor product of $d$ one-dimensional grids with levels $l_1,\ldots,l_d$.\footnote{Multiplications 
and inequalities have to be understood componentwise. We use the notation $\boldsymbol{1}=(1,\ldots,1)$.}
For each $\Omega_{\boldsymbol{l}}$, one considers hat functions $\phi_{\boldsymbol{l}, \boldsymbol{i}}$ centered at points $\boldsymbol{x}_{\boldsymbol{l},\boldsymbol{i}}$
defined by the product of the one-dimensional basis functions,
\begin{equation}
\phi_{\boldsymbol{l}, \boldsymbol{i}}(\boldsymbol{x}) = \prod_{j=1}^{d} \phi_{l_j, i_j}(x_j), \quad \boldsymbol{1} \leq \boldsymbol{i} \leq \boldsymbol{2^l-1}.
\label{phi}
\end{equation}

\noindent As in one dimension, one considers the function spaces spanned by these functions
\begin{equation}
V_{\boldsymbol{l}} = \mrm{span}\{\phi_{\boldsymbol{l},\boldsymbol{i}} : \boldsymbol{1}\leq \boldsymbol{i}\leq \boldsymbol{2^l-1}\}
\end{equation}

\noindent and the hierarchical increments
\begin{equation}
W_{\boldsymbol{l}} = \mrm{span}\{\phi_{\boldsymbol{l},\boldsymbol{i}} : \boldsymbol{i}\in \boldsymbol{I}_{\boldsymbol{l}}\},
\end{equation}

\noindent with $I_{\boldsymbol{l}} = \{\boldsymbol{i}\in\N^d : \boldsymbol{1}\leq \boldsymbol{i}\leq \boldsymbol{2^{l}-1}, \, i_j\,\text{odd for all $j$}\}$ and
\begin{equation}
V_{\boldsymbol{n}} = \bigoplus_{\boldsymbol{1}\leq \boldsymbol{l}\leq \boldsymbol{n}} W_{\boldsymbol{l}}.
\end{equation}

\noindent The \textit{multi-dimensional hierarchical basis} is the basis that corresponds to the $W_{\boldsymbol{l}}$'s for $\boldsymbol{1}\leq\boldsymbol{l}\leq\boldsymbol{n}$.
We show all subspaces $W_{\boldsymbol{l}}$ in two dimensions for $\boldsymbol{n}=(3,3)$ in Figure~\ref{sparsegrid}.

Equipped with a multi-dimensional hierarchical basis, one may approximate functions of $d$ variables. 
The appropriate function spaces in this context are the Korobov spaces $X^{2,p}(\Omega)$
defined for $2\leq p\leq+\infty$ by\footnote{For simplicity, we only consider functions that are zero at the boundary. Sparse grids for functions that are non-zero
at the boundary can be derived in an analogous fashion, and have similar approximation properties.}
\begin{equation}
X^{2,p}(\Omega) = \{f \in L^p(\Omega) : f\vert_{\partial\Omega}=0, \, D^{\boldsymbol{k}} f \in L^p(\Omega), \, \vert \boldsymbol{k}\vert_\infty \leq 2\},
\label{mixedSobolev}
\end{equation}

\noindent with $\vert \boldsymbol{k}\vert_\infty = \max_{1\leq j\leq d} k_j$ and norm
\begin{equation}
\vert f\vert_{\boldsymbol{2}, \infty} = \bigg\Vert\frac{\partial^{2d} f}{\partial x_1^{2} \ldots \partial x_d^{2}}\bigg\Vert_\infty.
\label{seminorm}
\end{equation}

\noindent These spaces go back to the 1959 paper of Korobov \cite{korobov1959}.
Note the difference with the Sobolev spaces $W^{2,p}(\Omega)$ defined in~\eqref{Sobolev}: 
smoothness for $X^{2,p}(\Omega)$ is measured in terms of \textit{mixed derivatives} of order two. 
For example in two dimensions, from $\vert\boldsymbol{k}\vert_\infty=\max(k_1,k_2)\leq2$,
one can see that the Korobov spaces $X^{2,p}(\Omega)$ require
\begin{equation}
\frac{\partial f}{\partial x_1},\,
\frac{\partial f}{\partial x_2},\,
\frac{\partial^2 f}{\partial x_1^2},\,
\frac{\partial^2 f}{\partial x_2^2},\,
\frac{\partial^2 f}{\partial x_1\partial x_2},\,
\frac{\partial^3 f}{\partial x_1^2\partial x_2},\,
\frac{\partial^3 f}{\partial x_1\partial x_2^2},\,
\frac{\partial^4 f}{\partial x_1^2\partial x_2^2} \in L^p(\Omega),
\end{equation}

\noindent whereas $\vert\boldsymbol{k}\vert_1=(k_1+k_2)\leq2$ for $W^{2,p}(\Omega)$ yields
\begin{equation}
\frac{\partial f}{\partial x_1},\,
\frac{\partial f}{\partial x_2},\,
\frac{\partial^2 f}{\partial x_1^2},\,
\frac{\partial^2 f}{\partial x_2^2},\,
\frac{\partial^2 f}{\partial x_1\partial x_2} \in L^p(\Omega).
\end{equation}

\noindent In other words, Korobov spaces $X^{2,p}(\Omega)$ are subsets of Sobolev spaces $W^{2,p}(\Omega)$.

The key fact is that any function $f\in X^{2,p}(\Omega)$ has a unique (infinite) expansion in the hierarchical basis,
\begin{equation}
f(\boldsymbol{x}) = \sum_{\boldsymbol{l}}\sum_{\boldsymbol{i}\in\boldsymbol{I_l}} v_{\boldsymbol{l}, \boldsymbol{i}}\phi_{\boldsymbol{l}, \boldsymbol{i}}(\boldsymbol{x}).
\label{expansion}
\end{equation}

\noindent The \textit{hierarchical coefficients} $v_{\boldsymbol{l}, \boldsymbol{i}}\in\R$ of $f$ are defined by~\cite[Lemma~3.2]{bungartz2004}
\begin{equation}
v_{\boldsymbol{l}, \boldsymbol{i}} = \int_\Omega \prod_{j=1}^d\Big(-2^{-(l_j+1)}\phi_{l_j,i_j}(x_j)\Big)\frac{\partial^{2d} f}{\partial x_1^{2} \ldots \partial x_d^{2}}(\boldsymbol{x})d\boldsymbol{x},
\end{equation}

\noindent and satisfy~\cite[Lemma~3.3]{bungartz2004}
\begin{equation}
\vert v_{\boldsymbol{l}, \boldsymbol{i}}\vert \leq 2^{-d}2^{-2\vert\boldsymbol{l}\vert_1}\vert f\vert_{\boldsymbol{2},\infty}.
\label{decaycoeffs}
\end{equation}

%%%%%%%%%%%%%%%%%%%%%%%%%%%%%%%%%%%%%%%%%%%%%%%%%%%%%%%%%%%%%%%%%%%%%%%%%%%
\subsection{Discretization}\label{section2.3}

\begin{figure}[t]
\centering
\begin{tikzpicture}
[   cnode/.style={draw=black, fill=#1, minimum width=3mm, circle},
]   
	% Dashed line:
	\draw[dashed] (7.75, 1.75) -- (7.75, -0.25) -- (5.75, -0.25) -- (5.75, -2.25) -- (3.75, -2.25) -- (3.75, -4.25) -- (1.75, -4.25) {};
	
	% Centers of grids (l2=1):
	\node[fill, circle, scale=0.2] at (2+1.5/2, 1.5/2) {};
	\node[fill, circle, scale=0.2] at (4+1.5/4, 1.5/2) {};
	\node[fill, circle, scale=0.2] at (4+4.5/4, 1.5/2) {};
 	\node[fill, circle, scale=0.2] at (6+1.5/8, 1.5/2) {};
	\node[fill, circle, scale=0.2] at (6+4.5/8, 1.5/2) {};
	\node[fill, circle, scale=0.2] at (6+7.5/8, 1.5/2) {};
	\node[fill, circle, scale=0.2] at (6+10.5/8, 1.5/2) {};
	
	% Centers of grid (l2=2):
	\node[fill, circle, scale=0.2] at (2+1.5/2, -1.25+1.5/4) {};
	\node[fill, circle, scale=0.2] at (2+1.5/2, -1.25-1.5/4) {};
	\node[fill, circle, scale=0.2] at (4+1.5/4, -1.25+1.5/4) {};
	\node[fill, circle, scale=0.2] at (4+4.5/4, -1.25+1.5/4) {};
	\node[fill, circle, scale=0.2] at (4+1.5/4, -1.25-1.5/4) {};
	\node[fill, circle, scale=0.2] at (4+4.5/4, -1.25-1.5/4) {};
	\node[fill, circle, scale=0.2] at (6+1.5/8, -1.25+1.5/4) {};
	\node[fill, circle, scale=0.2] at (6+4.5/8, -1.25+1.5/4) {};
	\node[fill, circle, scale=0.2] at (6+7.5/8, -1.25+1.5/4) {};
	\node[fill, circle, scale=0.2] at (6+10.5/8, -1.25+1.5/4) {};
	\node[fill, circle, scale=0.2] at (6+1.5/8, -1.25-1.5/4) {};
	\node[fill, circle, scale=0.2] at (6+4.5/8, -1.25-1.5/4) {};
	\node[fill, circle, scale=0.2] at (6+7.5/8, -1.25-1.5/4) {};
	\node[fill, circle, scale=0.2] at (6+10.5/8, -1.25-1.5/4) {};
	
	% Centers of grid (l2=3	):
	\node[fill, circle, scale=0.2] at (2+1.5/2, -3.25+4.5/8) {};
	\node[fill, circle, scale=0.2] at (2+1.5/2, -3.25+1.5/8) {};
	\node[fill, circle, scale=0.2] at (2+1.5/2, -3.25-1.5/8) {};
	\node[fill, circle, scale=0.2] at (2+1.5/2, -3.25-4.5/8) {};
	\node[fill, circle, scale=0.2] at (4+1.5/4, -3.25+4.5/8) {};
	\node[fill, circle, scale=0.2] at (4+1.5/4, -3.25+1.5/8) {};
	\node[fill, circle, scale=0.2] at (4+1.5/4, -3.25-1.5/8) {};
	\node[fill, circle, scale=0.2] at (4+1.5/4, -3.25-4.5/8) {};
	\node[fill, circle, scale=0.2] at (4+4.5/4, -3.25+4.5/8) {};
	\node[fill, circle, scale=0.2] at (4+4.5/4, -3.25+1.5/8) {};
	\node[fill, circle, scale=0.2] at (4+4.5/4, -3.25-1.5/8) {};
	\node[fill, circle, scale=0.2] at (4+4.5/4, -3.25-4.5/8) {};	
	\node[fill, circle, scale=0.2] at (6+1.5/8, -3.25+4.5/8) {};
	\node[fill, circle, scale=0.2] at (6+1.5/8, -3.25+1.5/8) {};
	\node[fill, circle, scale=0.2] at (6+1.5/8, -3.25-1.5/8) {};
	\node[fill, circle, scale=0.2] at (6+1.5/8, -3.25-4.5/8) {};
	\node[fill, circle, scale=0.2] at (6+4.5/8, -3.25+4.5/8) {};
	\node[fill, circle, scale=0.2] at (6+4.5/8, -3.25+1.5/8) {};
	\node[fill, circle, scale=0.2] at (6+4.5/8, -3.25-1.5/8) {};
	\node[fill, circle, scale=0.2] at (6+4.5/8, -3.25-4.5/8) {};
	\node[fill, circle, scale=0.2] at (6+7.5/8, -3.25+4.5/8) {};
	\node[fill, circle, scale=0.2] at (6+7.5/8, -3.25+1.5/8) {};
	\node[fill, circle, scale=0.2] at (6+7.5/8, -3.25-1.5/8) {};
	\node[fill, circle, scale=0.2] at (6+7.5/8, -3.25-4.5/8) {};
	\node[fill, circle, scale=0.2] at (6+10.5/8, -3.25+4.5/8) {};
	\node[fill, circle, scale=0.2] at (6+10.5/8, -3.25+1.5/8) {};
	\node[fill, circle, scale=0.2] at (6+10.5/8, -3.25-1.5/8) {};
	\node[fill, circle, scale=0.2] at (6+10.5/8, -3.25-4.5/8) {};
	
	% Legend:
	\node at (2.75, 2) {$l_1=1$};
	\node at (4.75, 2) {$l_1=2$};
	\node at (6.75, 2) {$l_1=3$};
	\node at (1.25, 0.75) {$l_2=1$};
	\node at (1.25, -1.25) {$l_2=2$};
	\node at (1.25, -3.25) {$l_2=3$};
	
	% l2= 1:
	\draw[draw=black] (2, 0) rectangle ++(1.5, 1.5);
	\draw[draw=black] (4, 0) rectangle ++(1.5/2, 1.5);
	\draw[draw=black] (4+1.5/2, 0) rectangle ++(1.5/2, 1.5);
	\draw[draw=black] (6, 0) rectangle ++(1.5/4, 1.5);
	\draw[draw=black] (6+1.5/4, 0) rectangle ++(1.5/4, 1.5);
	\draw[draw=black] (6+3/4, 0) rectangle ++(1.5/4, 1.5);
	\draw[draw=black] (6+4.5/4, 0) rectangle ++(1.5/4, 1.5);
	
 	% l2 = 2:
	\draw[draw=black] (2, -2) rectangle ++(1.5, 1.5/2);
	\draw[draw=black] (2, -2+1.5/2) rectangle ++(1.5, 1.5/2);
	\draw[draw=black] (4, -2) rectangle ++(1.5/2, 1.5/2);
	\draw[draw=black] (4+1.5/2, -2+1.5/2) rectangle ++(1.5/2, 1.5/2);
	\draw[draw=black] (4+1.5/2, -2) rectangle ++(1.5/2, 1.5/2);
	\draw[draw=black] (4, -2+1.5/2) rectangle ++(1.5/2, 1.5/2);
	\draw[draw=black] (6, -2) rectangle ++(1.5/4, 1.5/2);
	\draw[draw=black] (6+1.5/4, -2+1.5/2) rectangle ++(1.5/4, 1.5/2);
	\draw[draw=black] (6+1.5/4, -2) rectangle ++(1.5/4, 1.5/2);
	\draw[draw=black] (6, -2+1.5/2) rectangle ++(1.5/4, 1.5/2);
	\draw[draw=black] (6+3/4, -2) rectangle ++(1.5/4, 1.5/2);
	\draw[draw=black] (6+4.5/4, -2+1.5/2) rectangle ++(1.5/4, 1.5/2);
	\draw[draw=black] (6+4.5/4, -2) rectangle ++(1.5/4, 1.5/2);
	\draw[draw=black] (6+3/4, -2+1.5/2) rectangle ++(1.5/4, 1.5/2);
	
	% l2 = 3:
	\draw[draw=black] (2, -4) rectangle ++(1.5, 1.5);
	\draw[draw=black] (2, -4+1.5/4) rectangle ++(1.5, 1.5/4);
	\draw[draw=black] (2, -4+3/4) rectangle ++(1.5, 1.5/4);
	\draw[draw=black] (2, -4+4.5/4) rectangle ++(1.5, 1.5/4);
	\draw[draw=black] (4, -4) rectangle ++(1.5/2, 1.5);
	\draw[draw=black] (4, -4+1.5/4) rectangle ++(1.5/2, 1.5/4);
	\draw[draw=black] (4, -4+3/4) rectangle ++(1.5/2, 1.5/4);
	\draw[draw=black] (4, -4+4.5/4) rectangle ++(1.5/2, 1.5/4);
	\draw[draw=black] (4+1.5/2, -4) rectangle ++(1.5/2, 1.5);
	\draw[draw=black] (4+1.5/2, -4+1.5/4) rectangle ++(1.5/2, 1.5/4);
	\draw[draw=black] (4+1.5/2, -4+3/4) rectangle ++(1.5/2, 1.5/4);
	\draw[draw=black] (4+1.5/2, -4+4.5/4) rectangle ++(1.5/2, 1.5/4);
	\draw[draw=black] (6, -4) rectangle ++(1.5/4, 1.5/4);
	\draw[draw=black] (6+1.5/4, -4) rectangle ++(1.5/4, 1.5/4);
	\draw[draw=black] (6+3/4, -4) rectangle ++(1.5/4, 1.5/4);
	\draw[draw=black] (6+4.5/4, -4) rectangle ++(1.5/4, 1.5/4);
	\draw[draw=black] (6, -4+1.5/4) rectangle ++(1.5/4, 1.5/4);
	\draw[draw=black] (6+1.5/4, -4+1.5/4) rectangle ++(1.5/4, 1.5/4);
	\draw[draw=black] (6+3/4, -4+1.5/4) rectangle ++(1.5/4, 1.5/4);
	\draw[draw=black] (6+4.5/4, -4+1.5/4) rectangle ++(1.5/4, 1.5/4);
    \draw[draw=black] (6, -4+3/4) rectangle ++(1.5/4, 1.5/4);
	\draw[draw=black] (6+1.5/4, -4+3/4) rectangle ++(1.5/4, 1.5/4);
	\draw[draw=black] (6+3/4, -4+3/4) rectangle ++(1.5/4, 1.5/4);
	\draw[draw=black] (6+4.5/4, -4+3/4) rectangle ++(1.5/4, 1.5/4);
	\draw[draw=black] (6, -4+4.5/4) rectangle ++(1.5/4, 1.5/4);
	\draw[draw=black] (6+1.5/4, -4+4.5/4) rectangle ++(1.5/4, 1.5/4);
	\draw[draw=black] (6+3/4, -4+4.5/4) rectangle ++(1.5/4, 1.5/4);
	\draw[draw=black] (6+4.5/4, -4+4.5/4) rectangle ++(1.5/4, 1.5/4);
	
\end{tikzpicture}
\includegraphics[scale=0.33]{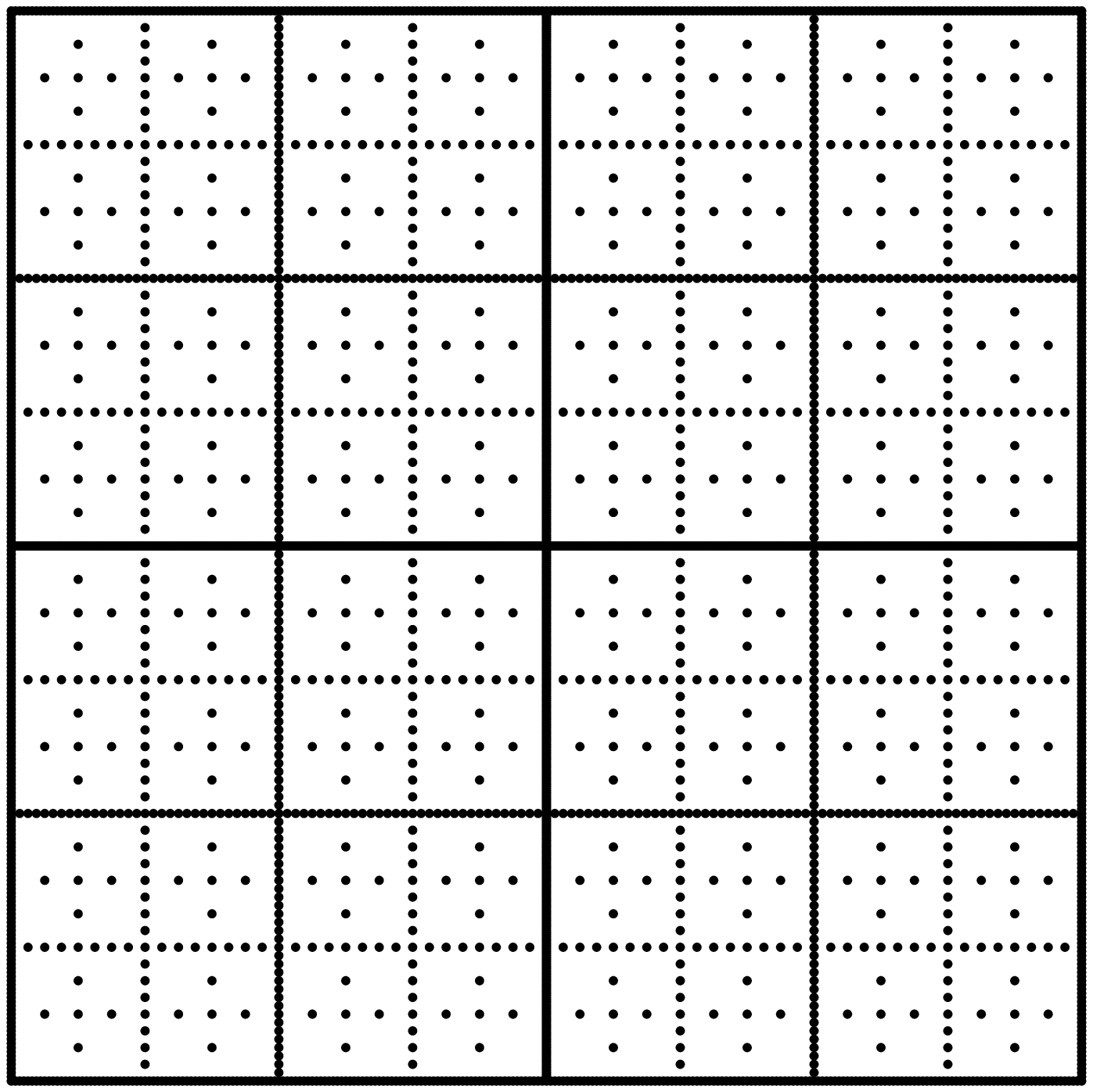}
\caption{\textit{Left: Subspaces $W_{\boldsymbol{l}}$ in two dimensions for $(l_1,l_2)\leq (3,3)$, and sparse and full grids $V_3^{(1)}$ and $V_3^{(\infty)}$.
The sparse grid consists of all subspaces above the dashed line $(\vert\boldsymbol{l}\vert_1=l_1+l_2\leq4)$.
The full grid consists of all subspaces $(\vert\boldsymbol{l}\vert_\infty=\max(l_1,l_2)\leq3)$.
By omitting many of the cross-terms that have small support in multiple dimensions, sparse grids have a significantly smaller number
of degrees of freedom than full grids, while their accuracy is only slightly deteriorated.
The dots represent the points at which the basis functions take the value $1$.
Each basis function lives inside a rectangle, i.e., for each $\boldsymbol{l}$ the supports of the hierarchical basis functions spanning 
$W_{\boldsymbol{l}}$ are mutually disjoint. Right: A sparse grid in two dimensions.}}
\label{sparsegrid}
\end{figure}

The third and last ingredient is a clever truncation of the expansion~\eqref{expansion}.
Sparse grids are discretizations of $X^{2,p}(\Omega)$ defined by
\begin{equation}
V_n^{(1)} = \bigoplus_{1\leq \vert\boldsymbol{l}\vert_1\leq n + d - 1} W_{\boldsymbol{l}},
\end{equation}

\noindent and correspond to a number of grid points $N$ given by~\cite[Lemma~3.6]{bungartz2004}
\begin{equation}
N =\mathcal{O}(h_n^{-1}\vert\log_2 h_n\vert^{d-1}) = \mathcal{O}(2^n n^{d-1}).
\end{equation}

\noindent A sparse grid in two dimensions is shown in Figure~\ref{sparsegrid}. 
Note that full grids $V_n^{(\infty)}$, with
\begin{equation}
V_n^{(\infty)} = \bigoplus_{1\leq \vert\boldsymbol{l}\vert_\infty\leq n} W_{\boldsymbol{l}},
\end{equation}

\noindent correspond to a much larger $\mathcal{O}(h_n^{-d}) = \mathcal{O}(2^{nd})$ number of grid points.

For any $f_n^{(1)}\in V_n^{(1)}$,
\begin{equation}
f_n^{(1)}(\boldsymbol{x}) = \sum_{\vert\boldsymbol{l}\vert_1\leq n+d-1}\sum_{\boldsymbol{i}\in\boldsymbol{I_l}} v_{\boldsymbol{l}, \boldsymbol{i}}\phi_{\boldsymbol{l}, \boldsymbol{i}}(\boldsymbol{x}),
\label{fn1}
\end{equation}

\noindent and for any $f\in X^{2,p}(\Omega)$, the approximation error satisfies~\cite[Lemma~3.13]{bungartz2004},
\begin{equation}
\Vert f - f^{(1)}_n\Vert_\infty = \mathcal{O}(N^{-2}\vert\log_2 N\vert^{3(d-1)}),
\label{approxerror}
\end{equation}

\noindent and, for any accuracy $\epsilon>0$,
\begin{equation}
\Vert f - f^{(1)}_n\Vert_\infty = \epsilon\quad\text{with}\quad N = \mathcal{O}(\epsilon^{-\frac{1}{2}}\vert\log_2\epsilon\vert^{\frac{3}{2}(d-1)}).
\label{decayerror}
\end{equation}

\noindent Th approximation error in \eqref{approxerror} is slightly worse than the $\mathcal{O}(N^{-\frac{2}{d}})$ error for approximating functions in $X^{2,p}(\Omega)$ 
with full grids~\cite[Lemma~3.5]{bungartz2004}, but using a much smaller number of points.

%%%%%%%%%%%%%%%%%%%%%%%%%%%%%%%%%%%%%%%%%%%%%%%%%%%%%%%%%%%%%%%%%%%%%%%%%%%
\section{Error bounds using sparse grids}\label{section3}

Results listed in Table~\ref{table:curse} are typically proven using the following technique. 
One shows that certain functions $f$ can be approximated by polynomials $f_M$ of degree $M$ to any
prescribed accuracy $\epsilon$, and so can polynomials $f_M$ by neural networks $f_N$ of size $N$, 
with $N$ bounded by some function of $\epsilon$, as in \eqref{approx}.
This amounts to decomposing the approximation error as 
\begin{equation}
\Vert f - f_N\Vert \leq \Vert f - f_M\Vert + \Vert f_M - f_N\Vert,
\end{equation}

\noindent for some norm $\Vert\cdot\Vert$. 
We use the same idea but instead of polynomials, we use approximations by sparse grids.

We first explain in Section \ref{section3.1} how to approximate the hat functions $\phi_{\boldsymbol{l},\boldsymbol{i}}$ using 
ideas introduced independently by Liang and Srikant \cite{liang2016}, and Yartosky \cite{yarotsky2017}.
We then prove in Section \ref{section3.2} our theorem concerning approximation of functions in $X^{2,p}(\Omega)$ by deep networks.

%%%%%%%%%%%%%%%%%%%%%%%%%%%%%%%%%%%%%%%%%%%%%%%%%%%%%%%%%%%%%%%%%%%%%%%%%%%
\subsection{Approximating multi-dimensional hat functions by deep networks}\label{section3.1}

The following proposition of Yarotsky shows how deep networks can implement multiplication.

\textbf{Proposition~1 \cite[Prop.~3]{yarotsky2017}.} \textit{For any $0<\epsilon<1$, there is a deep ReLU network with inputs $x_1$ and $x_2$,
with $\vert x_1\vert\leq M$ and $\vert x_2\vert\leq M$, that implements the multiplication $x_1x_2$ with accuracy $\epsilon$, outputs $0$ if $x_1=0$ or $x_2=0$,
and has depth and size $\mathcal{O}(\vert\log_2\epsilon\vert + \log_2 M)$.}

From Proposition~1, we obtain the following result, which shows how deep networks can approximate the multi-dimensional hat functions \eqref{phi} with a binary tree structure.
The corresponding network is shown in Figure~\ref{subnetwork}, and the proof can be found in Appendix A.

\textbf{Proposition~2.} \textit{For any dimension $d$ and $0<\epsilon<1$, there is a deep ReLU network with $d$ inputs $x_1,\ldots,x_d$ that implements the 
multiplication $\phi_{\boldsymbol{l},\boldsymbol{i}}(\boldsymbol{x})=\prod_{j=1}^d\phi_{l_j,i_j}(x_j)$ with accuracy $\epsilon$, outputs $0$ if one of the $\phi_{l_j,i_j}(x_j)$ is $0$, 
and has depth and size $\mathcal{O}(\vert\log_2\epsilon\vert\log_2 d)$.}

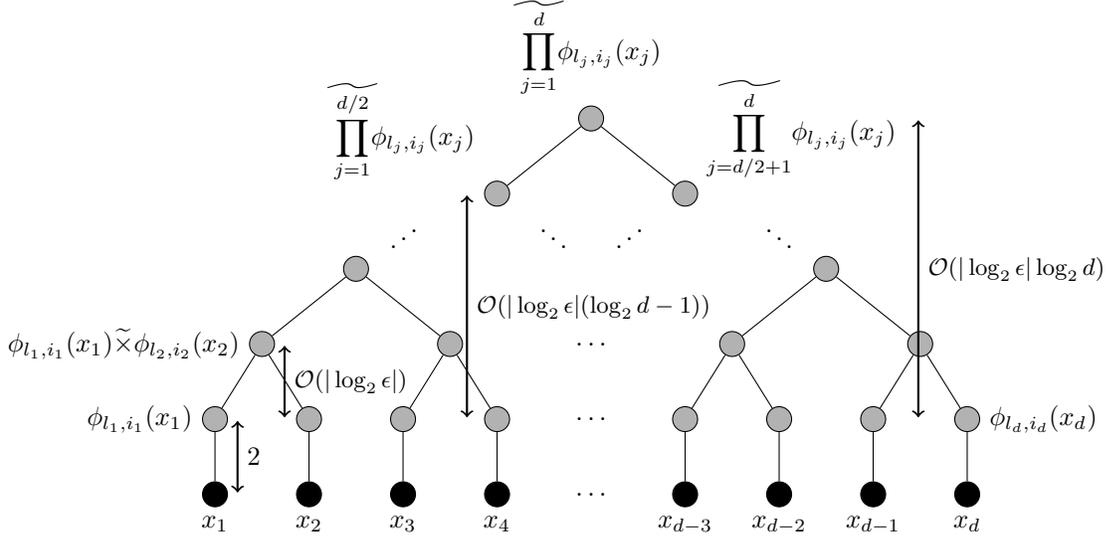
\begin{figure}
\centering
\begin{tikzpicture}
[   cnode/.style={draw=black, fill=#1, minimum width=3mm, circle},
]   
	% Inputs:
    \node[cnode = black, label = 270: $x_1$] (x1) at (1, -3) {};
    \node[cnode = black, label = 270: $x_2$] (x2) at (2.25, -3) {};
    \node[cnode = black, label = 270: $x_3$] (x3) at (3.5, -3) {};
    \node[cnode = black, label = 270: $x_4$] (x4) at (4.75, -3) {};
    \node at (6, -3) {$\ldots$};
    \node[cnode = black, label = 270: $x_{d-3}$] (xd3) at (7.25, -3) {};
    \node[cnode = black, label = 270: $x_{d-2}$] (xd2) at (8.5, -3) {};
    \node[cnode = black, label = 270: $x_{d-1}$] (xd1) at (9.75, -3) {};
    \node[cnode = black, label = 270: $x_{d}$] (xd) at (11, -3) {};
    
    % Phi functions:
    \node[cnode = gray, label = 180: $\phi_{l_1,i_1}(x_1)$] (phi1) at (1, -2) {};
    \node[cnode = gray] (phi2) at (2.25, -2) {};
    \node[cnode = gray] (phi3) at (3.5, -2) {};
    \node[cnode = gray] (phi4) at (4.75, -2) {};
    \node at (6, -2) {$\ldots$};
    \node[cnode = gray] (phid3) at (7.25, -2) {};
    \node[cnode = gray] (phid2) at (8.5, -2) {};
    \node[cnode = gray] (phid1) at (9.75, -2) {};
    \node[cnode = gray, label = 0: $\phi_{l_d,i_d}(x_d)$] (phid) at (11, -2) {};
    \draw (x1) -- (phi1) {};
    \draw (x2) -- (phi2) {};
    \draw (x3) -- (phi3) {};
    \draw (x4) -- (phi4) {};
    \draw (xd3) -- (phid3) {};
    \draw (xd2) -- (phid2) {};
    \draw (xd1) -- (phid1) {};
    \draw (xd) -- (phid) {};
    
    % Products:
    \node[cnode = gray, label = 180: $\phi_{l_1,i_1}(x_1)\widetilde{\times}\phi_{l_2,i_2}(x_2)$] (prod1) at (1.625, -1) {};
    \draw (phi1) -- (prod1) {};
    \draw (phi2) -- (prod1) {};
    \node[cnode = gray] (prod2) at (4.125, -1) {};
    \draw (phi3) -- (prod2) {};
    \draw (phi4) -- (prod2) {};
    \node at (6, -1) {$\ldots$};    
    \node[cnode = gray] (prod3) at (7.875, -1) {};
    \draw (phid3) -- (prod3) {};
    \draw (phid2) -- (prod3) {};
    \node[cnode = gray] (prod4) at (10.375, -1) {};
    \draw (phid1) -- (prod4) {};
    \draw (phid) -- (prod4) {};
    \node[cnode = gray] (prod5) at (2.875, 0) {};
    \draw (prod1) -- (prod5) {};
    \draw (prod2) -- (prod5) {};
    \node[cnode = gray] (prod6) at (9.125, 0) {};
    \draw (prod3) -- (prod6) {};
    \draw (prod4) -- (prod6) {};
    \node at (3.5, 0.5) {$\iddots$};
    \node at (5.5, 0.5) {$\ddots$};
    \node at (6.5, 0.5) {$\iddots$};
    \node at (8.5, 0.5) {$\ddots$};   
    \node[cnode = gray, label = 170: $\dsp\widetilde{\prod_{j=1}^{d/2}}\phi_{l_j,i_j}(x_j)$] (prod7) at (4.75, 1) {};
    \node[cnode = gray, label = 20: $\dsp\widetilde{\prod_{j=d/2+1}^{d}}\phi_{l_j,i_j}(x_j)$] (prod8) at (7.25, 1) {};
    \node[cnode = gray, label = 90: $\dsp\widetilde{\prod_{j=1}^{d}}\phi_{l_j,i_j}(x_j)$] (prod9) at (6, 2) {};
    \draw (prod7) -- (prod9) {};
    \draw (prod8) -- (prod9) {};
    
    % Depth of the network:
    \node[draw=none] (ghost1) at (1.925, -2.1) {};
    \node[draw=none] (ghost2) at (1.925, -0.9) {};
    \draw [thick, <->] (ghost1) -- node[right]{\small{$\mathcal{O}(\vert\log_2\epsilon\vert)$}} (ghost2);
    \node[draw=none] (ghost3) at (1.3, -3.1) {};
    \node[draw=none] (ghost4) at (1.3, -1.9) {};
    \draw [thick, <->] (ghost3) -- node[right]{$2$} (ghost4);    
    \node[draw=none] (ghost5) at (10.35, -2.1) {};
    \node[draw=none] (ghost6) at (10.35, 2.1) {};
    \draw [thick, <->] (ghost5) -- node[right]{\small{$\mathcal{O}(\vert\log_2\epsilon\vert\log_2d)$}} (ghost6);
    \node[draw=none] (ghost7) at (4.35, -2.1) {};
    \node[draw=none] (ghost8) at (4.35, 1.1) {};
    \draw [thick, <->] (ghost7) -- node[right]{\small{$\mathcal{O}(\vert\log_2\epsilon\vert(\log_2d-1))$}} (ghost8);
    
\end{tikzpicture}
\caption{\textit{A network that implements the $(d-1)$ products in $\prod_{j=1}^{d}\phi_{l_j,i_j}(x_j)$ with a binary tree structure.
Each product is computed with a subnetwork that has depth and size $\mathcal{O}(\vert\log_2\epsilon\vert)$.
The total product has accuracy $\epsilon$, and uses a network that has depth and size $\mathcal{O}(\vert\log_2\epsilon\vert\log_2d)$.}}
\label{subnetwork}
\end{figure}

%%%%%%%%%%%%%%%%%%%%%%%%%%%%%%%%%%%%%%%%%%%%%%%%%%%%%%%%%%%%%%%%%%%%%%%%%%%
\subsection{Approximating sparse grids by deep networks}\label{section3.2}

We use the fact that functions in $X^{2,p}([0,1]^d)$ can be approximated by sparse grids,
and then show that sparse grids can be represented by deep networks using the multiplication presented
in the previous subsection.
The resulting network is shown in Figure~\ref{sparsedeep}.

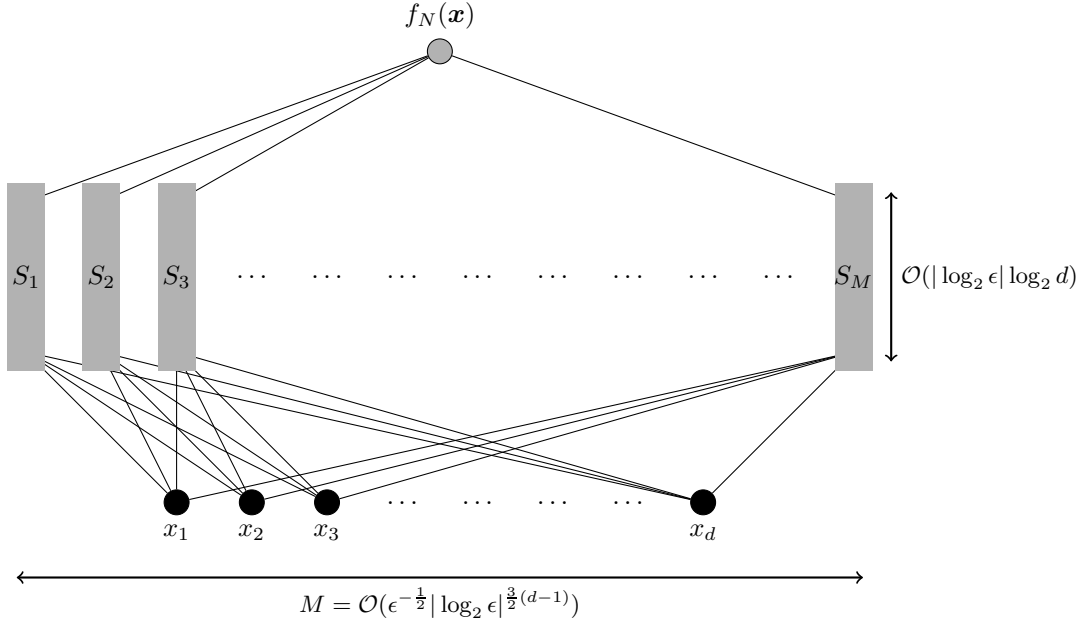
\begin{figure}
\centering
\begin{tikzpicture}
[   cnode/.style={draw=black, fill=#1, minimum width=3mm, circle},
]
	% Width of the network:
	\node[draw=none] (ghost1) at (-1.25, -4) {};
    \node[draw=none] (ghost2) at (10.25, -4) {};
    \draw [thick, <->] (ghost1) -- node[below]
    {\small{$M=\mathcal{O}(\epsilon^{-\frac{1}{2}}\vert\log_2\epsilon\vert^{\frac{3}{2}(d-1)})$}} (ghost2);
    
	% Inputs:
    \node[cnode = black, label = 270: $x_1$] (x1) at (1, -3) {};
    \node[cnode = black, label = 270: $x_2$] (x2) at (2, -3) {};
    \node[cnode = black, label = 270: $x_3$] (x3) at (3, -3) {};
    \node at (4, -3) {$\ldots$};
    \node at (5, -3) {$\ldots$};
    \node at (6, -3) {$\ldots$};
    \node at (7, -3) {$\ldots$};
    \node[cnode = black, label = 270: $x_d$] (xd) at (8, -3) {};
   
    % Subnetworks:
    \node[cnode = gray] (S1) at (-1, -1) {};
    \node[cnode = gray] (S2) at (0, -1) {};
    \node[cnode = gray] (S3) at (1, -1) {};
    \node[cnode = gray] (SN) at (10, -1) {};
    \draw (x1) -- (S1) {};
    \draw (x2) -- (S1) {};
    \draw (x3) -- (S1) {};
    \draw (xd) -- (S1) {};
    \draw (x1) -- (S2) {};
    \draw (x2) -- (S2) {};
    \draw (x3) -- (S2) {};
    \draw (xd) -- (S2) {};
    \draw (x1) -- (S3) {};
    \draw (x2) -- (S3) {};
    \draw (x3) -- (S3) {};
    \draw (xd) -- (S3) {};
    \draw (x1) -- (SN) {};
    \draw (x2) -- (SN) {};
    \draw (x3) -- (SN) {};
    \draw (xd) -- (SN) {};
    \node[cnode = gray] (SS1) at (-1,1) {};
    \node[cnode = gray] (SS2) at (0, 1) {};
    \node[cnode = gray] (SS3) at (1, 1) {};
    \node[cnode = gray] (SSN) at (10, 1) {};
    \node at (2, 0) {$\ldots$};
    \node at (3, 0) {$\ldots$};
    \node at (4, 0) {$\ldots$};
    \node at (5, 0) {$\ldots$};
    \node at (6, 0) {$\ldots$};
    \node at (7, 0) {$\ldots$};
    \node at (8, 0) {$\ldots$};
    \node at (9, 0) {$\ldots$};
    
    % Depth of the network:
    \node[draw=none] (ghost3) at (10.5, -1.25) {};
    \node[draw=none] (ghost4) at (10.5, 1.25) {};
    \draw [thick, <->] (ghost3) -- node[right]{\small{$\mathcal{O}(\vert\log_2\epsilon\vert\log_2d)$}} (ghost4);
    
    % Top of the network:
    \node[cnode = gray, label = 90: $f_N(\boldsymbol{x})$] (fN) at (4.5, 3) {};
    \draw (SS1) -- (fN) {};
    \draw (SS2) -- (fN) {};
    \draw (SS3) -- (fN) {};
    \draw (SSN) -- (fN) {};
    
    % Draw rectangles:
    \fill[gray] (-1.25,-1.25) -- (-0.75,-1.25) -- (-0.75,1.25) -- (-1.25,1.25) -- (-1.25,-1.25);
    \node at (-1, 0) {$S_1$};
    \fill[gray] (-0.25,-1.25) -- (0.25,-1.25) -- (0.25,1.25) -- (-0.25,1.25) -- (-0.25,-1.25);
    \node at (0, 0) {$S_2$};
    \fill[gray] (0.75,-1.25) -- (1.25,-1.25) -- (1.25,1.25) -- (0.75,1.25) -- (0.75,-1.25);
    \node at (1, 0) {$S_3$};
    \fill[gray] (9.75,-1.25) -- (10.25,-1.25) -- (10.25,1.25) -- (9.75,1.25) -- (9.75,-1.25);
    \node at (10, 0) {$S_M$};
    
\end{tikzpicture}
\caption{\textit{The sparse grid-based deep network used for the proof of Theorem~$1$. 
The network consists of $M$ subnetworks $S_1,S_2,\ldots,S_M$, which implement the multiplication $\prod_{j=1}^d\phi_{l_j,i_j}(x_j)$ 
presented in Section $\ref{section3.1}$.}}
\label{sparsedeep}
\end{figure}

\textbf{Theorem~1.} \textit{For any dimension $d$ and $0<\epsilon<1$, there is a deep ReLU network with $d$ inputs $x_1,\ldots,x_d$ 
capable of expressing any function $f$ in $X^{2,p}([0,1]^d)$ that satisfies $\vert f\vert_{\boldsymbol{2},\infty} \leq 1$ with accuracy $\epsilon$, 
and has depth $\mathcal{O}(\vert\log_2\epsilon\vert\log_2d)$ and size $\mathcal{O}(\epsilon^{-\frac{1}{2}}\vert\log_2\epsilon\vert^{\frac{3}{2}(d-1)+1}\log_2d)$.}

\textit{Proof.}
Let us consider $f\in X^{2,p}([0,1]^d)$ and suppose we want to approximate $f$ with a deep ReLU network $f_N$ of size $N$.
Let us write
\begin{equation}
\Vert f - f_N\Vert_\infty \leq \Vert f - f^{(1)}_m\Vert_\infty + \Vert f^{(1)}_m - f_N\Vert_\infty,
\end{equation}

\noindent where $f^{(1)}_m\in V_m^{(1)}$ is the sparse grid approximation of $f$ with $M = \mathcal{O}(2^mm^{d-1})$ points.
We know from \eqref{decayerror} that for any $\epsilon>0$, we can equal the first term to $\epsilon/2$ with
\begin{equation}
M = \mathcal{O}((\epsilon/2)^{-\frac{1}{2}}\vert\log_2\epsilon/2\vert^{\frac{3}{2}(d-1)}) 
= \mathcal{O}(\epsilon^{-\frac{1}{2}}\vert\log_2\epsilon\vert^{\frac{3}{2}(d-1)}).
\end{equation}

Let us now approximate $f_m^{(1)}$ by a network $f_N$ consisting of $M$ subnetworks, each subnetwork implementing the approximate
multiplication introduced in the previous subsection, which we write as $\widetilde{\phi}_{\boldsymbol{l}, \boldsymbol{i}}(\boldsymbol{x})$, that is,
\begin{equation}
f_N(\boldsymbol{x}) = \sum_{\vert\boldsymbol{l}\vert_1\leq m+d-1}\sum_{\boldsymbol{i}\in\boldsymbol{I_l}} v_{\boldsymbol{l}, \boldsymbol{i}}\widetilde{\phi}_{\boldsymbol{l}, \boldsymbol{i}}(\boldsymbol{x}).
\label{fN}
\end{equation}

\noindent Let us suppose that each $\widetilde{\phi}_{\boldsymbol{l}, \boldsymbol{i}}(\boldsymbol{x})$ is computed to accuracy $\delta$ with a network
of depth and size $\mathcal{O}(\vert\log_2\delta\vert\log_2d)$ for some $0<\delta<1$ (using Proposition~2).
From \eqref{fN}, we get
\begin{equation}
\vert f_m^{(1)}(\boldsymbol{x}) - f_N(\boldsymbol{x})\vert \leq \sum_{\vert\boldsymbol{l}\vert_1\leq m+d-1}\sum_{\boldsymbol{i}\in\boldsymbol{I_l}} \vert v_{\boldsymbol{l}, \boldsymbol{i}}\vert
\vert\phi_{\boldsymbol{l}, \boldsymbol{i}}(\boldsymbol{x})-\widetilde{\phi}_{\boldsymbol{l}, \boldsymbol{i}}(\boldsymbol{x})\vert.
\end{equation}

\noindent For a given $\boldsymbol{l}$, a given $\boldsymbol{x}$ belongs to the support of at most one $\phi_{\boldsymbol{l}, \boldsymbol{i}}(\boldsymbol{x})$ 
because these have disjoint supports,\footnote{Note that the same holds true for $\widetilde{\phi}_{\boldsymbol{l}, \boldsymbol{i}}(\boldsymbol{x})$ using the $0$-in-$0$-out property of Proposition 2.}
so the inequality becomes
\begin{equation}
\vert f_m^{(1)}(\boldsymbol{x}) - f_N(\boldsymbol{x})\vert \leq \sum_{\vert\boldsymbol{l}\vert_1\leq m+d-1}\vert v_{\boldsymbol{l}, \boldsymbol{i_l}}\vert
\vert\phi_{\boldsymbol{l}, \boldsymbol{i_l}}(\boldsymbol{x})-\widetilde{\phi}_{\boldsymbol{l}, \boldsymbol{i_l}}(\boldsymbol{x})\vert,
\end{equation}

\noindent for some $\boldsymbol{i_l}$, which yields, using $\vert\phi_{\boldsymbol{l}, \boldsymbol{i_l}}(\boldsymbol{x})-\tilde{\phi}_{\boldsymbol{l}, \boldsymbol{i_l}}(\boldsymbol{x})\vert\leq\delta$,
\begin{equation}
\vert f_m^{(1)}(\boldsymbol{x}) - f_N(\boldsymbol{x})\vert \leq \delta\sum_{\vert\boldsymbol{l}\vert_1\leq m+d-1}\vert v_{\boldsymbol{l}, \boldsymbol{i_l}}\vert.
\end{equation}

\noindent Using the property of the decay of the coefficients~\eqref{decaycoeffs}, 
$\sum_{\vert\boldsymbol{l}\vert_1\leq m+d-1}2^{-2\vert\boldsymbol{l}\vert_1}\leq 1$ and $2^{-d}\leq1$, we obtain
\begin{equation}
\vert f_m^{(1)}(\boldsymbol{x}) - f_N(\boldsymbol{x})\vert \leq\delta\vert f\vert_{\boldsymbol{2},\infty}\leq\delta,
\end{equation}

\noindent since $\vert f\vert_{\boldsymbol{2},\infty}\leq 1$. Hence, for $\delta = \epsilon/2$, one has 
\begin{equation}
\Vert f - f_N\Vert_\infty \leq \Vert f - f^{(1)}_m\Vert_\infty + \Vert f^{(1)}_m - f_N\Vert_\infty = \epsilon/2 + \epsilon/2 = \epsilon.
\end{equation}

\noindent The depth of the network is $\mathcal{O}(\vert\log_2\delta\vert\log_2d)=\mathcal{O}(\vert\log_2\epsilon\vert\log_2d)$, and its size is
\begin{equation}
N = \mathcal{O}(\vert\log_2\delta\vert\log_2d\times M) = \mathcal{O}(\epsilon^{-\frac{1}{2}}\vert\log_2\epsilon\vert^{\frac{3}{2}(d-1)+1}\log_2d).
\end{equation}

\noindent This completes the proof.\cqfd

%%%%%%%%%%%%%%%%%%%%%%%%%%%%%%%%%%%%%%%%%%%%%%%%%%%%%%%%%%%%%%%%%%%%%%%%%%%
\section{Discussion}

We have proven new rigorous upper bounds for the approximation of functions in Korobov spaces $X^{2,p}(\Omega)$ by deep ReLU networks, for which
the curse of dimensionality is lessened.
The proof is based on the ability of deep networks to approximate sparse grids via a binary tree structure (Figure \ref{subnetwork}), which resembles
the compositional structure used in \cite{poggio2017}.

There are many ways in which this work could be profitably continued.
To show an advantage of deep networks versus shallow, it would be desirable to obtain a lower bound for approximations in $X^{2,p}(\Omega)$ 
by shallow networks, for which the curse of dimensionality is \textit{not} lessened.\footnote{For approximations by shallow networks in Sobolev spaces $W^{2,2}(\Omega)$, 
it is well known that \textit{both} the lower and upper bounds depend exponentially on the dimension $d$ \cite[Th.~6.1]{pinkus1999}.}
Another extension would be to derive similar estimates for smoother functions, e.g., functions with mixed derivatives of order $m>2$. 
Piecewise smooth functions could also be considered (as in~\cite{petersen2017}), as well as Jacobi-weighted Korobov spaces \cite{shen2010}, 
and \textit{energy-based} sparse grids (for which the curse of dimensionality can be totally overcome~\cite[Th.~3.10]{bungartz2004}).
More generally, we could apply our methodology to any expansion of the form \eqref{expansion}, as long as 
the expansion coefficients satisfy a property like \eqref{decaycoeffs} (which controls the width of the network) and the
basis functions can be implemented efficiently using the multiplication of Section 3.1 (which controls the depth).

Our theorem provides an upper bound for the approximation complexity when the same network is used to approximate all functions
in a given Korobov space.
In other words, the network architecture does not depend on the function being approximated; only the weights $v_{\boldsymbol{l},\boldsymbol{i}}$ do. 
Alternatively, we could consider \textit{adaptive} architectures where not only the weights but also the architecture is adjusted to the function being approximated.
We would expect that this would decrease the complexity of the resulting network.
Adaptive network architectures in the context of approximating multivariate functions have been studied by, e.g., Yarotsky in \cite{yarotsky2017}.

As mentioned in the introduction, breaking the curse of dimensionality often relies on taking advantage of 
special properties of the functions being approximated. 
In this paper, we followed a different route and considered a more generic space of functions, and approximations by sparse grids.
Let us emphasize, however, that sparse grids---in particular the norm~\eqref{seminorm}---are highly anisotropic: 
to be efficient, these require the functions being approximated to be \textit{aligned with the axes}.
This is in fact the case for many algorithms for the approximation of multivariate functions,
including low-rank compressions and quasi-Monte Carlo methods; we refer the interested reader to~\cite{trefethen2017a} 
for details.

%%%%%%%%%%%%%%%%%%%%%%%%%%%%%%%%%%%%%%%%%%%%%%%%%%%%%%%%%%%%%%%%%%%%%%%%%%%
\section*{Acknowledgements}

We thank the members of the CM3 group (Ran Gu, Hwi Lee, Qi Sun and Yunzhe Tao) at Columbia University,
and colleague and programmer Joel R.~Clay for fruitful discussions.
The first author is much indebted to former PhD supervisor Nick Trefethen for his inspirational contributions
to numerical analysis and in particular to the field of approximation theory.

%%%%%%%%%%%%%%%%%%%%%%%%%%%%%%%%%%%%%%%%%%%%%%%%%%%%%%%%%%%%%%%%%%%%%%%%%%%
\bibliographystyle{siam}
\bibliography{/Users/montanelli/Dropbox/HM/WORK/ARTICLES_BOOKS/BIB_FILES/references.bib}

%%%%%%%%%%%%%%%%%%%%%%%%%%%%%%%%%%%%%%%%%%%%%%%%%%%%%%%%%%%%%%%%%%%%%%%%%%%
\appendix

\newpage
%%%%%%%%%%%%%%%%%%%%%%%%%%%%%%%%%%%%%%%%%%%%%%%%%%%%%%%%%%%%%%%%%%%%%%%%%%%
\section{Proof of Proposition 2}

Let us first note that $\phi_{l_j,i_j}$ can be written as
\begin{equation}
\phi_{l_j,i_j}(x_j) = \sigma\Bigg(1 - \sigma\bigg(\frac{x_j-x_{l_j,i_j}}{h_{l_j}}\bigg) - \sigma\bigg(\frac{x_{l_j,i_j}-x_j}{h_{l_j}}\bigg)\Bigg), 
\quad \sigma(x) = \max(0,x),
\end{equation}

\noindent i.e., it can be implemented by a network of depth $2$ and size $3$.

Let us now prove the result concerning the multiplication by induction over $d$. 
For simplicity we will suppose that $d=2^p$ is a power of $2$, and will prove the result by induction over $p$.

For $p=1$, i.e., $d=2$, we want to show that for any $0<\epsilon<1$ there is a deep ReLU network with two inputs $x_1$ and $x_2$ that implements the 
multiplication $\phi_{l_1,i_1}(x_1)\times\phi_{l_2,i_2}(x_2)$ with accuracy $\epsilon$, outputs $0$ if $\phi_{l_1,i_1}(x_1)=0$ or
$\phi_{l_2,i_2}(x_2)=0$ (which we call $0$-in-$0$-out property), and has depth and size $\mathcal{O}(\vert\log_2\epsilon\vert)$.
To create such a network, one can combine two networks that compute $\phi_{l_1,i_1}(x_1)$ and $\phi_{l_2,i_2}(x_1)$ from $x_1$ and $x_2$
(each having depth $2$ and size $3$) with the network of Proposition~1 (with $M=1$ since $\Vert \phi_{l_j,i_j}\Vert_\infty\leq 1$) to multiply $\phi_{l_1,i_1}(x_1)$ by $\phi_{l_2,i_2}(x_2)$ 
(depth and size $\mathcal{O}(\vert\log_2\epsilon\vert)$). 
The resulting network has depth $\mathcal{O}(2+\vert\log_2\epsilon\vert)=\mathcal{O}(\vert\log_2\epsilon\vert)$ 
and size $\mathcal{O}(3\times2+\vert\log_2\epsilon\vert)=\mathcal{O}(\vert\log_2\epsilon\vert)$, and inherits the $0$-in-$0$-out property from Proposition 1.

Let us suppose now that this is true in dimension $d/2=2^{p-1}$ for some $p\geq1$, and let us show this is still true in dimension $d=2^p$ for any $0<\epsilon<1$.
The induction hypothesis (which we use with $\epsilon/4$) states that there is a deep ReLU network with $d/2$ inputs $x_1,\ldots,x_{d/2}$ that implements
the multiplication $\prod_{j=1}^{d/2}\phi_{l_j,i_j}(x_j)$ with accuracy $\epsilon/4$, outputs $0$ if one of the $\phi_{l_j,i_j}(x_j)$ is $0$,
and has depth and size $\mathcal{O}(\vert\log_2\epsilon/4\vert(\log_2d-1)) = \mathcal{O}(\vert\log_2\epsilon\vert(\log_2d-1))$.
Let
\begin{equation}
\begin{array}{ll}
& \dsp \widetilde{\prod_{j=1}^{d/2}}\phi_{l_j,i_j}(x_j), \\[16pt]
= & \dsp \Bigg(\widetilde{\prod_{j=1}^{d/4}}\phi_{l_j,i_j}(x_j)\Bigg)\widetilde{\times}\Bigg(\widetilde{\prod_{j=d/4+1}^{d/2}}\phi_{l_j,i_j}(x_j)\Bigg), \\[16pt]
= & \dsp \Bigg(\widetilde{\prod_{j=1}^{d/8}}\phi_{l_j,i_j}(x_j)\widetilde{\times}\widetilde{\prod_{j=d/8+1}^{d/4}}\phi_{l_j,i_j}(x_j)\Bigg) 
\widetilde{\times}\Bigg(\widetilde{\prod_{j=d/4+1}^{3d/8}}\phi_{l_j,i_j}(x_j)\widetilde{\times}\widetilde{\prod_{j=3d/8+1}^{d/2}}\phi_{l_j,i_j}(x_j)\Bigg), \\[16pt]
= & \dsp \ldots,
\end{array}
\label{net1a}
\end{equation}

\noindent denote this network, where $\widetilde{\times}$ is the approximate multiplication of Proposition~1. 
In other words, this network corresponds to the hierarchical combination of $d/2-1$ products; see Figure \ref{subnetwork}.
The induction hypothesis tells us that this network has accuracy $\epsilon/4$, 
\begin{equation}
\Bigg\vert\widetilde{\prod_{j=1}^{d/2}}\phi_{l_j,i_j}(x_j) - \prod_{j=1}^{d/2} \phi_{l_j,i_j}(x_j)\Bigg\vert\leq\epsilon/4,
\label{ind1a}
\end{equation}

\noindent which yields
\begin{equation}
\Bigg\vert\widetilde{\prod_{j=1}^{d/2}}\phi_{l_j,i_j}(x_j)\Bigg\vert \leq \Bigg\vert\prod_{j=1}^{d/2} \phi_{l_j,i_j}(x_j)\Bigg\vert + \epsilon/4
\leq 1 + \epsilon/4.
\label{bound1a}
\end{equation}

\noindent Similarly we consider the network with $d/2$ inputs $x_{d/2+1},\ldots,x_{d}$ that implements $\prod_{j=d/2+1}^{d}\phi_{l_j,i_j}(x_j)$,
\begin{equation}
\widetilde{\prod_{j=d/2+1}^{d}}\phi_{l_j,i_j}(x_j) 
= \Bigg(\widetilde{\prod_{j=d/2+1}^{3d/4}}\phi_{l_j,i_j}(x_j)\Bigg)\widetilde{\times}\Bigg(\widetilde{\prod_{j=3d/4+1}^{d}}\phi_{l_j,i_j}(x_j)\Bigg)
= \ldots,
\label{net1b}
\end{equation}

\noindent with 
\begin{equation}
\Bigg\vert\widetilde{\prod_{j=d/2+1}^{d}}\phi_{l_j,i_j}(x_j) - \prod_{j=d/2+1}^{d} \phi_{l_j,i_j}(x_j)\Bigg\vert\leq\epsilon/4, \quad
\Bigg\vert\widetilde{\prod_{j=d/2+1}^{d}}\phi_{l_j,i_j}(x_j)\Bigg\vert \leq 1 + \epsilon/4,
\label{bound1b}
\end{equation}

\noindent and $0$-in-$0$-out property.
To construct a network that implements the full multiplication $\prod_{j=1}^{d}\phi_{l_j,i_j}(x_j)$, we combine \eqref{net1a} with \eqref{net1b}, that is,
\begin{equation}
\widetilde{\prod_{j=1}^{d}}\phi_{l_j,i_j}(x_j) = \Bigg(\widetilde{\prod_{j=1}^{d/2}}\phi_{l_j,i_j}(x_j)\Bigg)\widetilde{\times}\Bigg(\widetilde{\prod_{j=d/2+1}^{d}}\phi_{l_j,i_j}(x_j)\Bigg).
\label{net2}
\end{equation}

\noindent Note that \eqref{net2} satisfies the $0$-in-$0$-out property since \eqref{net1a} and \eqref{net1b} do.
Let us now examine the accuracy of this network:
\begin{equation}
\begin{array}{ll}
& \dsp \Bigg\vert\widetilde{\prod_{j=1}^{d}}\phi_{l_j,i_j}(x_j) - \prod_{j=1}^{d} \phi_{l_j,i_j}(x_j)\Bigg\vert, \\[16pt]
= & \dsp \Bigg\vert\Bigg(\widetilde{\prod_{j=1}^{d/2}}\phi_{l_j,i_j}(x_j)\Bigg)\widetilde{\times}\Bigg(\widetilde{\prod_{j=d/2+1}^{d}}\phi_{l_j,i_j}(x_j)\Bigg) 
- \Bigg(\prod_{j=1}^{d/2} \phi_{l_j,i_j}(x_j)\Bigg)\times\Bigg(\prod_{j=d/2+1}^{d} \phi_{l_j,i_j}(x_j)\Bigg)\Bigg\vert, \\[16pt]
\leq & \dsp\Bigg\vert\Bigg(\widetilde{\prod_{j=1}^{d/2}}\phi_{l_j,i_j}(x_j)\Bigg)\widetilde{\times}\Bigg(\widetilde{\prod_{j=d/2+1}^{d}}\phi_{l_j,i_j}(x_j)\Bigg)  - 
\Bigg(\widetilde{\prod_{j=1}^{d/2}}\phi_{l_j,i_j}(x_j)\Bigg)\times\Bigg(\widetilde{\prod_{j=d/2+1}^{d}}\phi_{l_j,i_j}(x_j)\Bigg) \Bigg\vert, \\[16pt]
+ & \dsp \Bigg\vert\Bigg(\widetilde{\prod_{j=1}^{d/2}}\phi_{l_j,i_j}(x_j)\Bigg)\times\Bigg(\widetilde{\prod_{j=d/2+1}^{d}}\phi_{l_j,i_j}(x_j)\Bigg)  
- \Bigg(\prod_{j=1}^{d/2} \phi_{l_j,i_j}(x_j)\Bigg)\times\Bigg(\prod_{j=d/2+1}^{d} \phi_{l_j,i_j}(x_j)\Bigg)\Bigg\vert.
\end{array}
\label{accuracy}
\end{equation}

\noindent Lets us first consider the second term:
\begin{equation}
\begin{array}{ll}
 & \dsp \Bigg\vert\Bigg(\widetilde{\prod_{j=1}^{d/2}}\phi_{l_j,i_j}(x_j)\Bigg)\times\Bigg(\widetilde{\prod_{j=d/2+1}^{d}}\phi_{l_j,i_j}(x_j)\Bigg)  
 - \Bigg(\prod_{j=1}^{d/2} \phi_{l_j,i_j}(x_j)\Bigg)\times\Bigg(\prod_{j=d/2+1}^{d} \phi_{l_j,i_j}(x_j)\Bigg)\Bigg\vert, \\[16pt]
\leq & \dsp\Bigg\vert\Bigg(\widetilde{\prod_{j=1}^{d/2}}\phi_{l_j,i_j}(x_j)\Bigg)\times\Bigg(\widetilde{\prod_{j=d/2+1}^{d}}\phi_{l_j,i_j}(x_j)\Bigg)  - 
\Bigg(\widetilde{\prod_{j=1}^{d/2}}\phi_{l_j,i_j}(x_j)\Bigg)\times\Bigg(\prod_{j=d/2+1}^{d}\phi_{l_j,i_j}(x_j)\Bigg) \Bigg\vert, \\[16pt]
+ & \dsp \Bigg\vert\Bigg(\widetilde{\prod_{j=1}^{d/2}}\phi_{l_j,i_j}(x_j)\Bigg)\times\Bigg(\prod_{j=d/2+1}^{d}\phi_{l_j,i_j}(x_j)\Bigg)  
 - \Bigg(\prod_{j=1}^{d/2} \phi_{l_j,i_j}(x_j)\Bigg)\times\Bigg(\prod_{j=d/2+1}^{d} \phi_{l_j,i_j}(x_j)\Bigg)\Bigg\vert, \\[16pt]
= & \dsp\Bigg\vert\widetilde{\prod_{j=1}^{d/2}}\phi_{l_j,i_j}(x_j)\Bigg\vert
\Bigg\vert\widetilde{\prod_{j=d/2+1}^{d}}\phi_{l_j,i_j}(x_j) - \prod_{j=d/2+1}^{d}\phi_{l_j,i_j}(x_j)\Bigg\vert, \\[16pt]
+ & \dsp \Bigg\vert\prod_{j=d/2+1}^{d} \phi_{l_j,i_j}(x_j)\Bigg\vert\Bigg\vert\widetilde{\prod_{j=1}^{d/2}}\phi_{l_j,i_j}(x_j) - \prod_{j=1}^{d/2} \phi_{l_j,i_j}(x_j)\Bigg\vert,\\[16pt]
\leq & (1+\epsilon/4)\times\epsilon/4 + 1\times\epsilon/4 = \epsilon/2 + \epsilon^2/16, \\
\end{array}
\end{equation}

\noindent using \eqref{ind1a}, \eqref{bound1a} and \eqref{bound1b}.
Therefore to bound \eqref{accuracy} by $\epsilon$ we would like to bound the first term in \eqref{accuracy} by $\epsilon/2 - \epsilon^2/16$.
Note that this term corresponds to the top multiplication of Figure \ref{subnetwork}.
To achieve accuracy $\epsilon/2 - \epsilon^2/16$, we use Proposition~1 with $M= 1 + \epsilon/4$.
Therefore, this multiplication is implemented by a network that has depth and size
\begin{equation}
\mathcal{O}(\vert(\log_2(\epsilon/2-\epsilon^2/16)\vert + \log_2(1+\epsilon/4)) 
= \mathcal{O}(\vert\log_2\epsilon\vert + 3\epsilon/8).
\end{equation}

\noindent The full network has depth and size
\begin{equation}
\mathcal{O}(\vert\log_2\epsilon\vert(\log_2d-1) + \vert\log_2\epsilon\vert +3\epsilon/8) 
= \mathcal{O}(\vert\log_2\epsilon\vert\log_2d).
\end{equation}

\noindent This completes the proof.\cqfd

%%%%%%%%%%%%%%%%%%%%%%%%%%%%%%%%%%%%%%%%%%%%%%%%%%%%%%%%%%%%%%%%%%%%%%%%%%%
\end{document}